\documentstyle[12pt,righttag,intlim]{amsart}
\title[Heat kernel asymptotics]
{Heat kernel asymptotics for Laplace type operators and matrix KdV hierarchy}
\author{Iosif Polterovich}\thanks{Supported by CRM-ISM postdoctoral fellowship}
\address{Institut des Sciences Math\'ematiques, Universit\'e du Qu\'ebec \`a
Montr\'eal\\ and  
Centre de Recherches Math\'ematiques, Universit\'e de Montr\'eal,
Montre\'al, Qu\'ebec, Canada}
\email{iossif@@math.uqam.ca}                                       
\dedicatory{{\rm Preliminary version}}
\def \phi{\varphi}
\def \epsilon{\varepsilon}
\numberwithin{equation}{subsection}
\theoremstyle{definition}
\newtheorem{definition}[equation]{Definition}
\theoremstyle{plain}
\newtheorem{lemma}[equation]{Lemma}
\newtheorem{theorem}[equation]{Theorem}

\begin{document}
\begin{abstract}
We study the heat kernel asymptotics for the Lap\-la\-ce type 
differential operators on vector bundles over Riemannian manifolds. 
In particular this includes the case of the  Laplacians acting on  
differential $p$-forms. We extend our results obtained earlier for
the scalar Laplacian and present closed formulas for all heat invariants 
associated with these operators. As another application, we present new
explicit formulas for the matrix Korteweg-de Vries hierarchy.
\end{abstract}
\maketitle
\section{Introduction and main results}
\subsection{Laplace type operators}
Let $M$ be a compact smooth $d$-dimensional Riemannian manifold
without boundary with a  metric $(g_{ij})$ and  $V$ be a smooth
$r$-dimensional vector bundle over $M$. Let $(x_1,\dots,x_d)$ be a
system of local coordinates centered at the point $x\in M$ and $(v_1,\dots,
v_r)$ be a local frame for $V$ defined near the origin. 

Consider a {\it Laplace type differential operator}, i.e. a 
second-order elliptic self-adjoint differential operator  
$D: C^{\infty}(V)\to C^{\infty}(V)$
with the scalar leading symbol given by the metric tensor (see [G1]). 
In the chosen local system it can be written as
\begin{equation}
\label{mat}
D=-(\sum_{i,j=1}^d g^{ij}I\cdot \partial^2/\partial x_i 
\partial x_j +\sum_{k=1}^d B_k \partial/\partial x_k + C),
\end{equation}
where $(g^{ij})$ denotes the inverse of the matrix $(g_{ij})$,
$I$ is the identity matrix and  $B_k=B_k(x_1,\dots,x_d)$, $C=C(x_1,\dots,x_d)$
are some endomorphisms of the  bundle $V$ given by $r\times r$ 
square matrices. 

In particular, any Laplacian
$\Delta^p$ acting on the  $\binom{d}{p}$-vector bundle $\Omega^p(M)$ 
of differential $p$-forms can be presented in this way (see [Ga]). 
Another example of a Laplace-type operator is a ($1$-dimensional) 
matrix Schr\"odinger operator (see section 3.1):
\begin{equation}
\label{schr}
L=I\cdot \frac{\partial^2}{\partial x^2} + U(x),
\end{equation}
where $U(x)$ is a $r\times r$ hermitian matrix.
\subsection{Heat kernel asymptotics}
Let us define the {\it heat operator} $e^{-tD}$ for $t>0$ 
(see [G1]). It is an infinitely smoothing operator from 
$L^2(V)\to C^{\infty}(V)$ and it is given by the kernel function 
$K(x,y,t, D): V_y \to V_x$ as follows:
$$(e^{-tD}f)(x)=\int_M K(x,y,t,D) f(y) \sqrt{g}dy.$$
where $g=\operatorname{det}(g_{ij})$.
The function $K(x,y,t,D)$ is called the {\it heat kernel} 
associated with the operator $D$. It has the following asymptotic
expansion on the diagonal as $t\to 0+$ (see [Se], [G1]):
$$K(x,x,t,D) \sim \sum_{n=0}^{\infty} A_n(x, D)t^{n-\frac{d}{2}},$$
where $A_n$ are certain endomorphisms of the fiber $V_x$.
  
Consider the spectral decomposition of the self-adjoing operator $D$
into a complete orthonormal basis of eigenfunctions $\phi_i$ and 
eigenvalues $\lambda_i$, $0 \le \lambda_1 \le \lambda_2 ...$.
Then
$$K(x,y,t,D)=\sum e^{-\lambda_i t} \phi_i(x)\otimes \phi_i(y)$$
Denoting $a_n(x, D)=\operatorname{Tr}(A_n(x, D))$ where $\operatorname{Tr}$
is the fiber trace of the endomorphism $A_n(x,D)$, we get
$$\operatorname{Tr}(K(x,x,t,D))=\sum_{i}e^{-t\lambda_i}
(\phi_i,\phi_i)(x) \sim \sum_{n=0}^{\infty} a_n(x, D)t^{n-\frac{d}{2}}.$$ 
The coefficients $a_n(x,D)$ are the {\it heat invariants} of the Riemannian
manifold $M$ associated with the operator $D$. 
They are homogeneous polynomials of degree $2n$ in the 
derivatives of the Riemannian metric $\{g^{ij}\}$ at the point $x$ ([G2]).

Computation of heat kernel coefficients is a long-standing problem
in spectral geometry (we refer to [P3] for a historical review).
In [P3] we calculated  heat invariants of Riemannian manifolds for the scalar 
Laplacian using the Agmon-Kannai (commutator) method (see [P1]). 
In this paper we present its generalization and extend our formulas for the 
heat invariants obtained in [P3] to the general Laplace type operators.
As another application, we extend the explicit expressions for  
the scalar Korteweg-de Vries hierarchy obtained in [P2], [P3] 
to the matrix case using the heat kernel asymptotics of the matrix
Schr\"odinger operator.

The main difficulty in applying the Agmon-Kannai method to matrix-valued 
operators is that the sum of the orders of two such operators is not 
necessarily greater than the order of their commutator (cf. (\ref{scal})).
However, for Laplace type operators this problem in fact does not appear (see
section 2.2).
\subsection{Main result}
Denote by  $\rho_x:M\to {\Bbb R}$  
the distance function on the manifold $M$: 
for every $y\in M$ the distance between the points $y$ and $x$ is 
$\rho_x(y)$. Let us also consider the operator $D$ and its powers as $r\times
r$ matrices (cf. (\ref{mat})). Let $\operatorname{Tr}(D^{j+n})$  
denote the scalar differential operator which is the matrix trace, 
i.e. the sum of diagonal entries of the operator $D^{j+n}$.
\begin{theorem}
\label{most}
The endomorphisms $A_n(x, D)$ are equal to 
\begin{equation*}
A_n(x, D)=(4\pi)^{-d/2}(-1)^n\sum_{j=0}^{3n}
\frac{\binom{3n+d/2}{j+d/2}}
{4^j \, j! \, (j+n)!}\left.
(D^{j+n})(\rho_x(y)^{2j}\cdot I)\right|_{y=x},
\end{equation*}
and hence the heat invariants $a_n(x, D)$ are given by
\begin{equation*}
a_n(x, D)=(4\pi)^{-d/2}(-1)^n\sum_{j=0}^{3n}
\frac{\binom{3n+d/2}{j+d/2}}
{4^j \, j! \, (j+n)!}\left.
\operatorname{Tr}(D^{j+n})(\rho_x(y)^{2j})\right|_{y=x}.
\end{equation*}
\end{theorem}
The binomial coefficients for $d$ odd are defined by
\begin{equation}
\binom{3n+d/2}{j+d/2}=
\frac{\Gamma(3n+1+d/2)}{(3n-j)! \, \Gamma(j+1+d/2)}.
\end{equation}
Theorem \ref{most} is proved in section 3.1.
\subsection{Remarks}
The endomorphisms $A_n(x,D)$ do not depend
on the choice of the coordinate system on $M$ and the local frame for $V$.
However, expressions in Theorem \ref{most} depend on
the frame in the fiber $V_x$. This means that these expressions
are subject to certain combinatorial cancellations which should reveal the 
invariant nature of the endomorphisms $A_n(x,D)$. Note that 
independence on the choice of the coordinate system is obtained exactly 
through the study of the combinatorial structure of our formulas (see [P3]).

Another indication of such hidden cancellations (as was noticed by P. Gilkey)
is the explicit dependence of the expressions in Theorem \ref{most} on the 
dimension $d$ of the manifold $M$. Indeed, if one writes $A_n(x,D)$ in  terms 
of  the curvature tensor and its derivatives, the dimension never appears.
Therefore, there should exist a way to rewrite Theorem \ref{most} in a 
dimension-independent form as well.
\section{The Agmon-Kannai method for Laplace type operators}
\subsection{The Agmon-Kannai expansion}
The method for computation of heat invariants developed in [P1], [P3] 
is based on an asymptotic expansion for resolvent kernels of
elliptic operators due to S.~Agmon and Y.~Kannai ([AK]). 
Let us briefly review the Agmon-Kannai result.
\begin{definition} (see [AK]).
Let $J=(j_1,\dots,j_l)$ be a  finite vector with nonnegative integer components
and let $P$, $Q$ be linear operators on a linear space $M$.
The {\it multiple commutator} $[P,Q;J]$  is defined for all
such vectors $J$ by
$$[P,Q;J]=
(\operatorname{ad} P)^{j_l}Q\cdots(\operatorname{ad} P)^{j_1}Q,$$
where $(\operatorname{ad} P)Q=[P,Q]$.
\end{definition} 
Let $H$ be a a self--adjoint elliptic differential 
operator of order $s$ on a manifold $M$ of dimension $d<s$, 
and $H_0$ be the operator obtained by
freezing the coefficients of the principal part $H'$ of the operator
$H$ at some point $x\in M$: $H_0=H'(x)$.  
Denote by $R_\lambda(x,y)$ the kernel of the resolvent
$R_\lambda=(H-\lambda)^{-1}$, and by $F_{\lambda}(x,y)$ --- the 
kernel of $F_{\lambda}=(H_0-\lambda)^{-1}$.
The Agmon--Kannai formula in its original form reads as follows ([AK]):
\begin{theorem}
\label{AK}
The following asymptotic representation on the diagonal 
holds for the kernel of the resolvent $R_{\lambda}$ 
as $\lambda \to \infty$:
\begin{equation}
\label{c}
R_\lambda(x,x) \sim \frac{1}{\sqrt{g}}(F_\lambda(x,x)+ 
\sum_{J}^\infty ([H_0-H,H_0;J]F^{|J|+r+1}(x,x))), 
\end{equation}
where the sum is taken over all vectors $J$ of length $\ge 1$ with 
nonnegative integer entries.
\end{theorem}
\subsection{An extension of the Agmon-Kannai theorem}
Denote by $o(P)$ the {\it order} of the differential operator $P$
(if $P=(P_{ij})$ is a matrix operator then $o(P)=\max_{i,j}o(P_{ij})$).
For any two scalar differential operators $P$, $Q$ we have:
\begin{equation}
\label{scal}
o([P,Q])\le o(P)+o(Q)-1
\end{equation}
This simple inequality lies in the basis of the Agmon-Kannai commutator 
expansions (cf. [AK], Lemma 5.1). Clearly, in the matrix case the relation 
(\ref{scal}) does not hold in general, and therefore Theorem \ref{AK} is
no longer true. However, under certain assumptions  Theorem \ref{AK} can be 
extended to the matrix case in a straightforward way (see [Kan] for an 
analogue of this theorem in a much more general situation).
\begin{lemma} 
\label{easy}
Let $H$ be an $r\times r$  matrix differential operator of order $s$ 
with a scalar principal part $H'$. Then the asymptotic expansion 
(\ref{c}) remains valid.
\end{lemma}
{\bf Proof.} 
Indeed, since $H'$ is a scalar operator, $H_0=h_0\cdot I$ is also scalar, 
$o(H_0)=s$.
Let $P$ be an $r\times r$ matrix differential of order $p$. 
Consider the  commutator $[P, H_0]$. Since $H_0$ is scalar, 
the commutator is a  matrix operator with the  entries  $[P_{ij}, h_0]$.  
For each entry we have: 
$o([P_{ij}, h_0]) \le p+s-1$, $i,j=1,\dots,r$. 
Therefore, $o([P, H_0]) \le p+s-1$. Let us note, that all commutators in
the formula (\ref{c}) are exactly of the  form $[P,H_0]$, where $P$ is
some matrix operator. We have shown that inequality (\ref{scal}) holds
for all commutators of this kind, and therefore in this case 
we may just repeat the proof of Theorem \ref{AK} (see [AK]).  
This completes the proof of the lemma. \qed

In [P1] we found a concise reformulation of the Agmon--Kannai theorem
for scalar operators (Theorem 1.2). Lemma \ref{easy} allows to extend it 
immediately to the case of matrix operators with a scalar principal 
part.
\begin{theorem}
\label{new}
Let $H$ be an elliptic differential operator of order $s$ acting on
a $r$-dimensional vector bundle $V$ over a Riemannian
manifold $M$ of order $d < s$.
The following asymptotic representation on the diagonal holds
for the resolvent kernel $R_\lambda(x,y)$ of an elliptic operator $H$ 
as $\lambda \to \infty$ :
\begin{equation}
\label{s}
R_\lambda(x,x) \sim \frac{1}{\sqrt{g}} 
\sum_{m=0}^\infty X_m F_\lambda^{m+1}(x,x),
\end{equation}
where the operators $X_m$ are defined by:
$$
X_m=\sum_{k=0}^m (-1)^k \binom{m}{k} H^k H_0^{m-k}.
$$
\end{theorem}
Note that the condition $d<s$ can be relaxed in the same way as in [P3] --- 
by taking the derivatives of the resolvent with respect to the spectral     
parameter $\lambda$ (see Theorem 2.3.1 in [P3]).
\subsection{Proof of Theorem 1.3.1} 
Indeed, consider a normal coordinate system 
centered at the point $x \in M$. Then $g_{ij}(x)=\delta_{ij}$ and hence
the principal part of the operator $D$ with the coefficients frozen at the
point $x$ has the form $D_0=\Delta_0 \cdot I$, where $\Delta_0$ is the scalar 
Laplacian on $R^d$ and $I$ is the $r\times r$ identity matrix.
Due to Theorem \ref{new} we can repeat the arguments of the proofs of 
Theorems 1.2.1 and 3.1.1 from [P3] just taking into account multiplication
by the matrix $I$. Finally we obtain the expression for $A_n(x, D)$
identical to the formula (4.2.2) in [P3]. 
Expressing the distance function in terms of normal coordinates as in [P3]
we obtain the formulas for the endomorphisms $A_n(x,D)$.
Taking the traces of $A_n(x,D)$ we get the heat invariants $a_n(x, D)$ 
and this completes the proof of the theorem.
\qed
\section{Matrix KdV hierarchy}
\subsection{Heat kernel asymptotics for Schr\"odinger operator} 
Consider the $1$-dimensional matrix Schr\"odinger operator $L$ 
(see (\ref{schr}))
with a potential $U$ which is an hermitian $r\times r$ matrix. 
Its heat kernel $H(t,x,y)$ is  
the fundamental solution of the heat equation 
$$\left(\frac{\partial}{\partial t} - L\right)f=0.$$
It has the following  asymptotic representation on the diagonal 
as $t\to 0+$: 
\begin{equation*}
H(t,x,x)\sim\frac{1}{\sqrt{4\pi t}}\sum_{n=0}^{\infty}h_n[U]t^n,
\end{equation*}
where $h_n[U]$ are some polynomials in the matrix $U(x)$ and 
its derivatives. 

The matrix KdV hierarchy is defined by (see [AvSc1]):
\begin{equation}
\label{kdvh}
\frac{\partial U}{\partial t}=\frac{\partial}{\partial x}G_n[U],
\end{equation}
where
\begin{equation*}
G_n[U]=
\frac{(2n)!}{2\cdot n!}h_n[U], 
\quad n\in {\Bbb N}.
\end{equation*}
Set $U_0=U$, $U_n=\partial^n U/\partial x^n$, $n\in {\Bbb N}$,
where $U_n$, $n\ge 0$ are formal variables. 
The sequence of polynomials $G_n[U]=G_n[U_0,U_1,U_2,\dots]$ 
starts with (see [AvSc2]): 
\begin{multline*}
G_1[U]=U_0, \quad
G_2[U]=U_2+3U_0^2, \\ 
G_3[u]=U_4+5U_0U_2+5U_2U_0+5U_1^2+10U_0^3, \, \dots
\end{multline*}

\subsection{Computation of matrix KdV hierarchy}
In [P2], [P3] we presented explicit 
formulas for the scalar Korteweg-de Vries hierarchy using heat invariants
of the scalar $1$-dimensional  Schr\"odinger operator
(we refer to [P2] for the history of this question). 
Theorem \ref{most} allows to extend our results to the matrix KdV hierarchy 
as well.

For other approaches to explicit computations of the matrix KdV hierarchy see 
[AvSc1], [AvSc2].
\begin{theorem}
\label{check}
The KdV hierarchy is  given by:
$$
G_n[U]=\frac{(2n)!}{2\cdot n!}\sum_{j=0}^{n}
\binom{n+\frac{1}{2}}{j+\frac{1}{2}}
\frac{(-1)^j}{4^j \, j! \, (j+n)!}P_{nj}[U],
$$
where the polynomial $P_{nj}[U]$ is obtained from 
$\left.L^{j+n}(x^{2j})\right|_{x=0}$ by a formal change of variables:
$U_i(0)\to U_i$, $i=0,..,2n+2j-2$.
\end{theorem}

This expression can be completely expanded due to a formula 
for the powers of the Schr\"odinger operator ([Rid]) which remains
valid in the matrix case as well (note that now $U_i$ are matrices and hence
do not commute).

\begin{theorem}
\label{KdV1}
The  polynomials $G_n[U]$, $n\in {\Bbb N}$ are equal to:
\begin{equation*}
G_n[U]=
\frac{(2n)!}{2\cdot n!}\sum_{j=0}^{n}
\binom{n+\frac{1}{2}}{j+\frac{1}{2}}
\frac{(-1)^j (2j)!}{4^j \, j! \, (j+n)!}
\sum_{p=1}^{j+n}\mskip-2\thinmuskip
\sum_{k_1,\dots,k_p\atop k_1+\cdots +k_p=2(n-p)}\mskip-13\thinmuskip
C_{k_1,\dots,k_p}U_{k_1}\cdots U_{k_p},
\end{equation*}
where 
\begin{equation*}
C_{k_1,\dots,k_p}=
\sum\begin{Sb}
0\le\l_0\le l_1\le \cdots\le l_{p-1}=j+n-p\\ 
2l_i\ge k_1+\cdots+k_{i+1},\,\,i=0,\dots,p-1.\end{Sb}
\binom{2l_0}{k_1}\binom{2l_1-k_1}{k_2}\cdots 
\binom{2l_{p-1}-k_1-\cdots-k_{p-1}}{k_p}.
\end{equation*}
\end{theorem}      
\subsection*{Acknowledgments}
I am very grateful to my Ph.D. advisor Yakar Kannai for  suggesting the 
subject of this research and constant support.
I would also like to thank Peter Gilkey and Leonid Polterovich for useful
discussions.
\section*{References}
\noindent [AK] S. Agmon, Y. Kannai, On the asymptotic behavior of spectral 
functions and resolvent kernels of elliptic operators, Israel J. Math. 5
(1967), 1-30.





\noindent [AvSc1] I.G. Avramidi, R. Schimming, A new explicit expression
for the Korteweg - de Vries hierarchy, solv-int/9710009,  (1997), 1-17.

\noindent [AvSc2] I.G. Avramidi, R. Schimming, Heat kernel coefficients
for the matrix Schr\"odinger operator, J. Math. Phys., 36 (1995), no. 9,
5042-5054.

\noindent [Ga] M.P. Gaffney, Asymptotic distributions associated with the
Laplacian on forms, Comm. Pure Appl. Math., vol. 11 (1958), 535-545.
 









\noindent [G1] P. Gilkey, The spectral geometry of a Riemannian manifold, J. 
Diff. Geom. 10 (1975), 601-618.

\noindent [G2]  P. Gilkey, The index theorem and the heat equation, Math. 
Lect.Series,  Publish or Perish, 1974.

\noindent [G3] P. Gilkey, Invariance theory, the heat equation and the Atiyah-
Singer index theorem, Mathematics Lecture Series, 11, Publish or Perish, 1984.

\noindent 
[Kan] Y. Kannai, On the asymptotic behavior of resolvent kernels, spectral
functions and eigenvalues of semi-elliptic systems, Ann. Scuola Norm. Sup.
Pisa (3), 23 (1969), 563-634.














\noindent [P1] I. Polterovich, A commutator method for computation of heat
invariants,  Indag. Math., N.S., (2000) 11(1) 139-149.


\noindent [P2] I. Polterovich, From Agmon-Kannai expansion to Korteweg-de Vries
hierarchy, Lett. Math. Phys. 49 (1999), 71-77.


\noindent [P3] I. Polterovich, Heat invariants of Riemannian manifolds,
to appear in Israel Journal of Mathematics.

\noindent [Rid] S.Z. Rida, Explicit formulae for the powers of a 
Schr\"odinger-like
ordinary differential operator, J. Phys. A: Math. Gen. 31 (1998), 5577-5583.




\noindent 
[Se] R. Seeley, Complex powers of an elliptic operator, Proc. Symp. Pure
Math. 10 (1967), 288-307.






\end{document}